

\baselineskip=14pt
\parskip=10pt
\def\halmos{\hbox{\vrule height0.15cm width0.01cm\vbox{\hrule height
  0.01cm width0.2cm \vskip0.15cm \hrule height 0.01cm width0.2cm}\vrule
  height0.15cm width 0.01cm}}

\magnification=\magstephalf

\def\1{{\overline{1}}}
\def\2{{\overline{2}}}
\parindent=0pt
\overfullrule=0in

\def\frac#1#2{{#1 \over #2}}
\centerline
{\bf Two Quick Proofs of a Catalan Lemma Needed by Lisa Sauermann and Yuval Wigderson }
\bigskip
\centerline
{\it Shalosh B. EKHAD and Doron ZEILBERGER}

\bigskip

 \qquad {\it In memory of Robin Chapman (1963-2020), a problem-solving maestro and a great Catalanist}

In a recent beautiful article [SW], the authors needed the following identity (Claim 3.8 there)

$$
\sum_{i=0}^{s} (-1)^i \, C_i {{i+1} \choose {s-i}} \, = \, 0  \quad ,
\eqno(1)
$$
where $C_i$ are the Catalan numbers.  Let's give two shorter proofs of $(1)$.

{\bf Proof 1} (by SBE): Go to Maple, and type [Note that one does {\bf not} need Zeilberger, Gosper suffices].

{\tt sum((-1)**i*binomial(2*i,i)/(i+1)*binomial(i+1,s-i),i=0..s); } \halmos

{\bf Proof 2} (by DZ) The Catalan number $C_i$ famously (inter alia), counts the number of (complete) binary trees
with $i+1$ leaves. ${{i+1} \choose {s-i}}$ counts the number of words with $s-i$ twos and  $2i+1-s$  ones, hence
the total number of words in the alphabet $\{1,2\}$ with $i+1$ letters whose sum is $s+1$.

Hence $\sum_{i=0}^{s}  \, C_i {{i+1} \choose {s-i}} $ counts {\it all} binary trees with $\leq s+1$ leaves
where the leaves have labels $1$ or $2$ that add-up to $s+1$. The lemma is equivalent to the fact
that the number of such creatues with an {\bf even} number of leaves 
equals the number of those with an {\bf odd} number of leaves. The following {\it involution} provides the needed bijection.

{\it Scan the leaves from left to right until you either encounter a leaf labeled $2$, in which case you make it
into an internal vertex and make it give birth to two new leaves each labeled $1$, or you encounter a leaf labeled $1$ whose
sister is also labeled $1$, in which case you remove them both and make their mother a new leaf labeled $2$.} \halmos

The reason they needed $(1)$ was to prove the following rather hairy identity
$$
\sum_{
{ 
{m_1, \dots, m_t  \geq 1}
\atop
{m_1 \, +   \, \dots \,+ \, m_t=l}
} 
}
\, (-1)^t {{l-m_1} \choose {m_1-1}} {{l-m_2} \choose {m_2}} \cdots {{l-m_t} \choose {m_t}} \, = \,
(-1)^l C_{l-1} \quad .
\eqno(2)
$$

Rather than using $(1)$ we will present two direct, shorter proofs. In fact we will prove the more general identity, for $l \geq m$,
$$
\sum_{
{ 
{m_1, \dots, m_t  \geq 1}
\atop
{m_1 \, +   \, \dots \,+ \, m_t=m}
} 
}
\, (-1)^t {{l-m_1} \choose {m_1-1}} {{l-m_2} \choose {m_2}} \cdots {{l-m_t} \choose {m_t}} \, = \,
(-1)^m C_{m-1} \quad .
\eqno(2')
$$
Identity $(2)$ is the special case $l=m$ of $(2')$.

Calling the left side of $(2')$ $A(l,m)$, separating the case $t=1$ (that yields $-{{l-m} \choose {m-1}}$)
and summing over $m_t$ (let's call it $k$), we readily get the recurrence
$$
A(l,m)= -{{l-m} \choose {m-1}} - \sum_{k=1}^{m-1} {{l-k} \choose {k}}\, A(l,m-k) \quad .
$$
It would then follow by induction on $m$ that $A(l,m)$ equals $(-1)^m C_{m-1}$, if the latter satisfies
the same recurrence. But this is equivalent to the  identity
$$
\sum_{k=0}^{m} (-1)^k \, {{l-m+k} \choose {m-k}} \,  C_k = \, {{l-m-1} \choose {m}}\quad .
\eqno(3)
$$

{\bf Proof 1} (by SBE): Dividing by the right hand side, this is, in turn, equivalent to
$$
\sum_{k=0}^{m} (-1)^k \, \frac{ {{l-m+k} \choose {m-k}}}{ {{l-m-1} \choose {m}}} \, C_k\, = \,1 \quad .
$$
Calling the left side $f(m)$, go into Maple and type

{\tt
SumTools[Hypergeometric][ZeilbergerRecurrence]((-1)**k*binomial(2*k,k)/(k+1)*binomial(l-m+k,m-k)/binomial(l-m-1,m),m,k,f,0..m);}

and in one nano-second you would get that $f(m)$ satisfies the following recurrence
$$
- \left( m+1 \right) f \left( m \right) + \left( m+2 \right) f \left( m+1 \right) =1 \quad,
$$
Since $f(1)=1$ (check!), it follows by induction that $f(m)=1$ for every $m$ (and of course every $l \geq m$). \halmos

{\bf Proof 2} (by DZ): The left side of $(3)$ without the $(-1)^k$ is the number of pairs 
$(T,w)$ where $T$ is a binary tree with $\leq m+1$ leaves and $w$ is a word in $\{1,2\}$ of length $l-m-1$  longer than
the number of leaves of $T$, whose sum is $l$. Hence the left side of $(3)$ is the difference between the
number of such creatures where $T$ has an odd number of leaves and those that have an even number of leaves.
The above involution that proved $(1)$ (where the leaves of $T$ are labeled by the corresponding prefix of $w$)
is still valid, but now the survivors are the pairs $( ., 1 w)$ where $.$ is
the one-leaf tree, and w is a word in $\{1,2\}$ of length $l-m-1$ whose sum is $l-1$. \halmos

{\bf Thanks} are due to Victor S. Miller for bringing [SW] to out attention. Also thanks to
Lisa Sauermann for very insightful and useful comments on an earlier draft.

\vfill\eject

{\bf Reference}

[SW] Lisa Sauermann and Yuval Wigderson, {\it Polynomials that vanish to high order on most of the hypercube},
arXiv:2010.00077v1 [math. C0] 30 Sep 2020

\bigskip
\hrule
\bigskip
Shalosh B. Ekhad and Doron Zeilberger, Department of Mathematics, Rutgers University (New Brunswick), Hill Center-Busch Campus, 110 Frelinghuysen
Rd., Piscataway, NJ 08854-8019, USA. \hfill\break
Email: {\tt [ShaloshBEkhad, DoronZeil] at gmail dot com}   \quad .

{\bf Exclusively published in the Personal Journal of Shalosh B. Ekhad and Doron Zeilberger and arxiv.org}

First written:{\bf Nov. 11, 2020} . This version: {\bf Nov. 15, 2020}.

\end